\documentclass[a4paper,12pt,twoside]{article}
\usepackage[francais]{babel}
\usepackage{amssymb, amsmath, eucal}
\usepackage[all]{xy}
\usepackage{mathrsfs}
\usepackage[dvips]{graphics}
\begin{document}
\textwidth15.5cm
\textheight22.5cm
\voffset=-13mm
\newtheorem{The}{Theorem}[section]
\newtheorem{Lem}[The]{Lemma}
\newtheorem{Prop}[The]{Proposition}
\newtheorem{Cor}[The]{Corollary}
\newtheorem{Rem}[The]{Remark}
\newtheorem{Obs}[The]{Observation}
\newtheorem{SConj}[The]{Standard Conjecture}
\newtheorem{Titre}[The]{\!\!\!\! }
\newtheorem{Conj}[The]{Conjecture}
\newtheorem{Question}[The]{Question}
\newtheorem{Prob}[The]{Problem}
\newtheorem{Def}[The]{Definition}
\newtheorem{Not}[The]{Notation}
\newtheorem{Claim}[The]{Claim}
\newtheorem{Conc}[The]{Conclusion}
\newtheorem{Ex}[The]{Example}
\newtheorem{Fact}[The]{Fact}
\newcommand{\C}{\mathbb{C}}
\newcommand{\R}{\mathbb{R}}
\newcommand{\N}{\mathbb{N}}
\newcommand{\Z}{\mathbb{Z}}
\newcommand{\Q}{\mathbb{Q}}
\newcommand{\Proj}{\mathbb{P}}

\begin{center}

{\Large\bf Sufficient Bigness Criterion for Differences of Two Nef Classes}

\end{center}

\begin{center}

{\large Dan Popovici}

\end{center}

\vspace{1ex}

\noindent {\small {\bf Abstract.} We prove the qualitative part of Demailly's conjecture on transcendental Morse inequalities for differences of two nef classes satisfying a numerical relative positivity condition on an arbitrary compact K\"ahler (and even more general) manifold. The result improves on an earlier one by J. Xiao whose constant $4n$ featuring in the hypothesis is now replaced by the optimal and natural $n$. Our method follows arguments by Chiose as subsequently used by Xiao up to the point where we introduce a new way of handling the estimates in a certain Monge-Amp\`ere equation. This result is needed to extend to the K\"ahler case and to transcendental classes the Boucksom-Demailly-Paun-Peternell cone duality theorem if one is to follow these authors' method and was conjectured by them.}

\section{Introduction}

 Let $X$ be a compact complex manifold with $\mbox{dim}_{\C}X=n$. Following various authors (e.g. [Chi13]), Xiao makes in [Xia13] the following assumption\!\!: \\

$(H)$\hspace{2ex} {\it there exists a Hermitian metric $\omega$ on $X$  such that} \\

\hspace{15ex}  {\it  $\partial\bar\partial\omega^k=0$ for all $k=1, 2, \dots , n-1$.} \\

\noindent It is clear that $(H)$ holds if $X$ is a K\"ahler manifold. It is also standard and easy to check that condition $(H)$ is equivalent to either of the following two equivalent conditions\!\!:

$$\partial\bar\partial\omega=0 \hspace{1ex} \mbox{and} \hspace{1ex} \partial\bar\partial\omega^2=0 \iff \partial\bar\partial\omega=0 \hspace{1ex} \mbox{and} \hspace{1ex} \partial\omega\wedge\bar\partial\omega=0.$$

\noindent Following Xiao's method in [Xia13], itself inspired by earlier authors, especially Chiose [Chi13], we prove the following statement in which a real Bott-Chern class of bidegree $(1,\,1)$ being {\it nef} means, as usual, that it contains $C^{\infty}$ representatives with arbitrarily small negative parts (see inequalities (\ref{eqn:varphi-epsilon-pos})).

\begin{The}\label{The:main} Let $X$ be a compact complex manifold with $\mbox{dim}_{\C}X=n$ satisfying the assumption $(H)$. Then, for any {\bf nef} Bott-Chern cohomology classes $\{\alpha\},\,\{\beta\}\in H^{1,\,1}_{BC}(X,\,\R)$, the following implication holds\!\!:

$$\{\alpha\}^n - n\,\{\alpha\}^{n-1}.\,\{\beta\} >0  \Longrightarrow \mbox{the class}\,\,\{\alpha - \beta\}\,\, \mbox{contains a K\"ahler current}.$$

\end{The}

 This answers affirmatively the {\it qualitative} part of a special version (i.e. the one for a difference of two nef classes) of Demailly's {\it transcendental Morse inequalities conjecture} (see [BDPP13, Conjecture 10.1, $(ii)$]) and will be crucial to the eventual extension of the duality theorem proved in [BDPP13, Theorem 2.2] to transcendental classes in the fairly general context of compact K\"ahler (not necessarily projective) manifolds. Although the method we propose here also produces a lower bound for the volume of the difference class $\{\alpha - \beta\}$, this bound (that we will not present here) is weaker than the lower bound $\{\alpha\}^n - n\,\{\alpha\}^{n-1}.\,\{\beta\}$ predicted in the {\it quantitative} part of Conjecture 10.1, $(ii)$ in [BDPP13]. 
 
 Xiao proves in [Xia13] the existence of a K\"ahler current in the class $\{\alpha - \beta\}$ under the stronger assumption $\{\alpha\}^n - 4n\,\{\alpha\}^{n-1}.\,\{\beta\} >0$ and the same assumption $(H)$ on $X$. The two ingredients he uses are as follows.

\begin{Lem}(Lamari's duality lemma, [Lam99, Lemme 3.3])\label{Lem:Lamari} Let $X$ be a compact complex manifold with $\mbox{dim}_{\C}X=n$ and let $\alpha$ be any $C^{\infty}$ real $(1,\,1)$-form on $X$. Then the following two statements are equivalent. \\

\noindent $(i)$\, There exists a distribution $\psi$ on $X$ such that $\alpha + i\partial\bar\partial\psi\geq 0$ in the sense of $(1,\,1)$-currents on $X.$ \\

\noindent $(ii)$\, $\displaystyle\int\limits_X\alpha\wedge\gamma^{n-1}\geq 0$ for any Gauduchon metric $\gamma$ on $X.$

\end{Lem}

As an aside, we notice that this statement, when applied to $d$-closed real $(1,\,1)$-forms $\alpha$, translates to the pseudo-effective cone ${\cal E}_X\subset H^{1,\,1}_{BC}(X,\,\R)$ of $X$ and the closure of the Gauduchon cone $\overline{{\cal G}_X}\subset H^{n-1,\,n-1}_A(X,\,\R)$ of $X$ being dual under the duality between the Bott-Chern cohomology of bidegree $(1,\,1)$ and the Aeppli cohomology of bidegree $(n-1,\,n-1)$. (See [Pop13] for the definition of the Gauduchon cone.)

\begin{The}(the Tosatti-Weinkove resolution of Hermitian Monge-Amp\`ere equations, [TW10]) \label{The:TW}  Let $X$ be a compact complex manifold with $\mbox{dim}_{\C}X=n$ and let $\omega$ be a Hermitian metric on $X$. 

 Then, for any $C^{\infty}$ function $F\,:\,X\rightarrow \R$, there exist a unique constant $C>0$ and a unique $C^{\infty}$ function $\varphi\,:\,X\rightarrow \R$ such that

$$(\omega + i\partial\bar\partial\varphi)^n = C e^F\omega^n,\hspace{2ex} \omega + i\partial\bar\partial\varphi >0 \hspace{2ex} \mbox{and} \hspace{2ex} \sup\limits_X\varphi=0.$$

\end{The}

 As a matter of fact, Yau's classical theorem that solved the Calabi Conjecture, of which Theorem \ref{The:TW} is a generalisation to the possibly non-K\"ahler context, suffices for the proof of Theorem \ref{The:main} whose assumptions imply that $X$ must be K\"ahler (as already pointed out by Xiao in his situation based on [Chi13, Theorem 0.2]) although this is not used either here or in Xiao's work.

\section{Xiao's approach}\label{section:Xiao} In this section, we simply reproduce Xiao's arguments (themselves inspired by earlier authors) up to the point where we will branch off in a different direction in the next section to handle certain estimates. 

 Let us fix a Hermitian metric $\omega$ on $X$ such that $\partial\bar\partial\omega^k=0$ for all $k$. We also fix nef Bott-Chern $(1,\,1)$-classes $\{\alpha\}, \{\beta\}$. By the nef assumption, for every $\varepsilon>0$, there exist $C^{\infty}$ functions $\varphi_{\varepsilon}, \psi_{\varepsilon}\,:\,X\rightarrow\R$ such that

\begin{equation}\label{eqn:varphi-epsilon-pos}\alpha_{\varepsilon}:=\alpha + \varepsilon\,\omega + i\partial\bar\partial\varphi_{\varepsilon}>0  \hspace{2ex}\mbox{and}\hspace{2ex}   \beta_{\varepsilon}:=\beta + \varepsilon\,\omega + i\partial\bar\partial\psi_{\varepsilon}>0 \hspace{2ex}\mbox{on}\hspace{2ex} X.\end{equation}

\noindent Note that $\alpha_{\varepsilon}$ and $\beta_{\varepsilon}$ need not be $d$-closed, but the property $\partial\bar\partial\omega^k=0$ yields\!\!:

\begin{equation}\label{eqn:varphi-epsilon-k-ddbar}\partial\bar\partial\alpha_{\varepsilon}^k = \partial\bar\partial\beta_{\varepsilon}^k=0 \hspace{2ex}\mbox{and}\hspace{2ex}  \partial\bar\partial(\alpha + \varepsilon\,\omega)^k = \partial\bar\partial(\beta + \varepsilon\,\omega)^k=0\end{equation}

\noindent for all $k=1, 2, \dots , n-1$. We normalise $\sup\limits_X\varphi_{\varepsilon} = \sup\limits_X\psi_{\varepsilon}=0$ for every $\varepsilon>0$.

 Let us fix $\varepsilon>0$. The existence of a K\"ahler current in the class $\{\alpha - \beta\} = \{\alpha_{\varepsilon} - \beta_{\varepsilon}\}$ is equivalent to

$$\exists\,\delta>0,\,\exists\,\,\mbox{a distribution}\,\theta_{\delta}\,\,\mbox{on}\,X \hspace{2ex} \mbox{such that} \hspace{2ex} \alpha_{\varepsilon} - \beta_{\varepsilon} + i\partial\bar\partial\theta_{\delta}\geq\delta\,\alpha_{\varepsilon},$$

\noindent which, in view of Lamari's duality lemma \ref{Lem:Lamari}, is equivalent to

$$\exists\delta>0 \hspace{2ex} \mbox{such that} \hspace{2ex} \int\limits_X(\alpha_{\varepsilon} - \beta_{\varepsilon})\wedge\gamma^{n-1}\geq\delta\,\int\limits_X\alpha_{\varepsilon}\wedge\gamma^{n-1}$$

\noindent for every Gauduchon metric $\gamma$ on $X$. This is, of course, equivalent to

$$\exists\delta>0 \hspace{2ex} \mbox{such that} \hspace{2ex} (1-\delta)\,\int\limits_X\alpha_{\varepsilon}\wedge\gamma^{n-1}\geq\int\limits_X\beta_{\varepsilon}\wedge\gamma^{n-1}$$

\noindent for every Gauduchon metric $\gamma$ on $X$.

 Xiao's approach is to prove the existence of a K\"ahler current in the class $\{\alpha - \beta\} = \{\alpha_{\varepsilon} - \beta_{\varepsilon}\}$ by contradiction. Suppose that no such current exists. Then, for every $\varepsilon>0$ and every sequence of positive reals $\delta_m\downarrow 0$, there exist Gauduchon metrics $\gamma_{m,\,\varepsilon}$ on $X$ such that

\begin{equation}\label{eqn:contrad-ineq} (1-\delta_m)\,\int\limits_X\alpha_{\varepsilon}\wedge\gamma_{m,\,\varepsilon}^{n-1} < \int\limits_X\beta_{\varepsilon}\wedge\gamma_{m,\,\varepsilon}^{n-1} =1 \hspace{3ex} \mbox{for all} \hspace{1ex} m\in\N^{\star},\,\varepsilon>0.\end{equation}

\noindent The last identity is a normalisation of the Gauduchon metrics $\gamma_{m,\,\varepsilon}$ which is clearly always possible by rescaling $\gamma_{m,\,\varepsilon}$ by a positive factor. This normalisation implies that for every $\varepsilon>0$, the positive definite $(n-1,\,n-1)$-forms $(\gamma_{m,\,\varepsilon}^{n-1})_m$ are uniformly bounded in mass, hence after possibly extracting a subsequence we can assume the convergence $\gamma_{m,\,\varepsilon}^{n-1}\rightarrow\Gamma_{\infty,\,\varepsilon}$ in the weak topology of currents as $m\rightarrow +\infty$, where $\Gamma_{\infty,\,\varepsilon}\geq 0$ is an $(n-1,\,n-1)$-current on $X$. Taking limits as $m\rightarrow +\infty$ in (\ref{eqn:contrad-ineq}), we get

\begin{equation}\label{eqn:alpha-epsilon-int-bound}\int\limits_X\alpha_{\varepsilon}\wedge\Gamma_{\infty,\,\varepsilon}\leq 1  \hspace{3ex} \mbox{for all} \hspace{1ex} \varepsilon>0.\end{equation}

 Note that the l.h.s. of (\ref{eqn:contrad-ineq}) does not change if $\alpha_{\varepsilon}$ is replaced with any $\alpha_{\varepsilon} + i\partial\bar\partial u$ (thanks to $\gamma_{m,\,\varepsilon}$ being Gauduchon), while $\alpha_{\varepsilon}\wedge\gamma_{m,\,\varepsilon}^{n-1}$ is (after division by $\gamma_{m,\,\varepsilon}^n$) the trace of $\alpha_{\varepsilon}$ w.r.t. $\gamma_{m,\,\varepsilon}$ divided by $n$ (i.e. the arithmetic mean of the eigenvalues). To find a lower bound for the trace that would contradict (\ref{eqn:contrad-ineq}), it is natural to prescribe the volume form (i.e. the product of the eigenvalues) of some $\alpha_{\varepsilon} + i\partial\bar\partial u_{m,\,\varepsilon}$ by imposing that it be, up to a constant factor, the strictly positive $(n,\,n)$-form featuring in the r.h.s. of (\ref{eqn:contrad-ineq}). More precisely, the Tosatti-Weinkove theorem \ref{The:TW} allows us to solve the Monge-Amp\`ere equation

$$(\star)_{m,\,\varepsilon}  \hspace{6ex} (\alpha_{\varepsilon} + i\partial\bar\partial u_{m,\,\varepsilon})^n = c_{\varepsilon}\,\beta_{\varepsilon}\wedge\gamma_{m,\,\varepsilon}^{n-1}$$

\noindent for any $\varepsilon>0$ and any $m\in\N^{\star}$ by ensuring the existence of a unique constant $c_{\varepsilon}>0$ and of a unique $C^{\infty}$ function $u_{m,\,\varepsilon}\,:\,X\rightarrow\R$ satisfying $(\star)_{m,\,\varepsilon}$ so that

$$\widetilde{\alpha}_{m,\,\varepsilon}:=\alpha_{\varepsilon} + i\partial\bar\partial u_{m,\,\varepsilon} >0, \hspace{2ex} \sup\limits_X(\varphi_{\varepsilon} + u_{m,\,\varepsilon}) = 0.$$

 Note that $c_{\varepsilon}$ is independent of $m$ since we must have

\begin{equation}\label{eqn:c_0-def}c_{\varepsilon} = \int\limits_X\widetilde{\alpha}_{m,\,\varepsilon}^n = \int\limits_X(\alpha + \varepsilon\omega)^n\downarrow\int\limits_X\alpha^n:=c_0>0,\end{equation}

\noindent where the non-increasing convergence is relative to $\varepsilon\downarrow 0$. Indeed, the second identity in (\ref{eqn:c_0-def}) follows from $\partial\bar\partial(\alpha + \varepsilon\omega)^k=0$ for all $k=1, 2, \dots , n-1$ (cf. (\ref{eqn:varphi-epsilon-k-ddbar})). Thus, it is significant that $c_{\varepsilon}$ does not change if we add any $i\partial\bar\partial u$ to $\alpha$, i.e. $c_{\varepsilon}$ depends only on the Bott-Chern class $\{\alpha\}$, on $\omega$ and on $\varepsilon$. Analogously, one defines

\begin{equation}\label{eqn:M_0-def}M_{\varepsilon}: = \int\limits_X\widetilde{\alpha}_{m,\,\varepsilon}^{n-1}\wedge\beta_{\varepsilon} = \int\limits_X(\alpha + \varepsilon\omega)^{n-1}\wedge(\beta + \varepsilon\omega)\downarrow\int\limits_X\alpha^{n-1}\wedge\beta:=M_0\geq 0,\end{equation}

\noindent where the non-increasing convergence is relative to $\varepsilon\downarrow 0$. Clearly, $M_{\varepsilon}$ is independent of $m$ and depends only on the Bott-Chern classes $\{\alpha\}$, $\{\beta\}$, on $\omega$ and on $\varepsilon$. Note that the second integral in (\ref{eqn:M_0-def}) equals $\int_X(\alpha + \varepsilon\omega + i\partial\bar\partial\varphi_{\varepsilon})^{n-1}\wedge(\beta + \varepsilon\omega + i\partial\bar\partial\psi_{\varepsilon})$ which is positive since $\alpha_{\varepsilon}, \beta_{\varepsilon}>0$ by (\ref{eqn:varphi-epsilon-pos}). Since $M_0\geq 0$, the hypothesis $c_0 - nM_0>0$ made in Theorem \ref{The:main} implies $c_0>0$. This justifies the final claim in (\ref{eqn:c_0-def}).

\section{Estimates in the Monge-Amp\`ere equation}\label{section:estimates}

 We now propose an approach to the details of these estimates that differs from that of Xiao. We start with a very simple, elementary (and probably known) observation.

\begin{Lem}\label{Lem:trace-prod-ineq} For any Hermitian metrics $\alpha, \beta, \gamma$ on a complex manifold, the following inequality holds at every point\!\!:

\begin{equation}\label{eqn:trace-prod-ineq}(\Lambda_{\alpha}\beta)\cdot(\Lambda_{\beta}\gamma) \geq \Lambda_{\alpha}\gamma.\end{equation}

\end{Lem}

\noindent {\it Proof.} Since (\ref{eqn:trace-prod-ineq}) is a pointwise inequality, we fix an arbitrary point $x$ and choose local coordinates about $x$ such that

\vspace{1ex}

\noindent $\beta(x) = \sum\limits_j i dz_j\wedge d\bar{z}_j, \hspace{2ex} \alpha(x) = \sum\limits_j \alpha_j\,i dz_j\wedge d\bar{z}_j \hspace{2ex} \mbox{and} \hspace{2ex} \gamma(x) = \sum\limits_{j,\,k} \gamma_{j\bar{k}}\,i dz_j\wedge d\bar{z}_k.$

\vspace{1ex}

\noindent Then $\alpha_j>0$ and $\gamma_{j\bar{j}}>0$ for every $j$. If we denote by the same symbol any $(1,\,1)$-form and its coefficient matrix in the chosen coordinates, we have

\vspace{1ex}

\hspace{10ex}  $\alpha^{-1}\,\gamma = (\frac{1}{\alpha_j}\,\gamma_{j\bar{k}})_{j,\,k}$, hence $\mbox{Tr}(\alpha^{-1}\,\gamma) = \sum\limits_j\frac{1}{\alpha_j}\,\gamma_{j\bar{j}}.$

\vspace{1ex}

\noindent Thus (\ref{eqn:trace-prod-ineq}) translates to $(\sum\limits_j\frac{1}{\alpha_j})\,\sum\limits_k\gamma_{k\bar{k}} \geq \sum\limits_j\frac{1}{\alpha_j}\,\gamma_{j\bar{j}}$ which clearly holds since $\sum\limits_{j\neq k}\frac{1}{\alpha_j}\,\gamma_{k\bar{k}} >0$ because all the $\alpha_j$ and all the $\gamma_{k\bar{k}}$ are positive.  \hfill $\Box$

\vspace{2ex}

Our main observation is the following statement.

\begin{Lem}\label{Lem:crucial-estimate} For every $m\in\N^{\star}$ and every $\varepsilon>0$, we have\!\!:

\begin{equation}\label{eqn:crucial-estimate}\bigg(\int\limits_X\widetilde{\alpha}_{m,\,\varepsilon}\wedge\gamma_{m,\,\varepsilon}^{n-1}\bigg)\cdot\bigg(\int\limits_X\widetilde{\alpha}_{m,\,\varepsilon}^{n-1}\wedge\beta_{\varepsilon}\bigg) \geq \frac{1}{n}\,\int\limits_X\widetilde{\alpha}_{m,\,\varepsilon}^n = \frac{c_{\varepsilon}}{n}.\end{equation}

\end{Lem}

\noindent {\it Proof.} Let $0<\lambda_1\leq\lambda_2\leq\dots\leq\lambda_n$, resp. $0<\mu_1\leq\mu_2\leq\dots\leq\mu_n$, be the eigenvalues of $\widetilde{\alpha}_{m,\,\varepsilon}$, resp. $\beta_{\varepsilon}$, w.r.t. $\gamma_{m,\,\varepsilon}$. We have\!\!: \\

\noindent $\displaystyle\widetilde{\alpha}_{m,\,\varepsilon}^n = \lambda_1\dots\lambda_n\,\gamma_{m,\,\varepsilon}^n  \hspace{1ex} \mbox{and} \hspace{1ex} \widetilde{\alpha}_{m,\,\varepsilon}\wedge\gamma_{m,\,\varepsilon}^{n-1} = \frac{1}{n}\,(\Lambda_{\gamma_{m,\,\varepsilon}}\widetilde{\alpha}_{m,\,\varepsilon})\,\gamma_{m,\,\varepsilon}^n = \frac{\lambda_1 + \dots + \lambda_n}{n}\, \gamma_{m,\,\varepsilon}^n.$ \\

\noindent Similarly, $\displaystyle\beta_{\varepsilon}\wedge\gamma_{m,\,\varepsilon}^{n-1} = \frac{1}{n}\,(\Lambda_{\gamma_{m,\,\varepsilon}}\beta_{\varepsilon})\,\gamma_{m,\,\varepsilon}^n = \frac{\mu_1 + \dots + \mu_n}{n}\, \gamma_{m,\,\varepsilon}^n.$ \\

\noindent Thus, the Monge-Amp\`ere equation $(\star)_{m,\,\varepsilon}$ translates to

\begin{equation}\label{eqn:eqn-translation}\lambda_1\dots\lambda_n = c_{\varepsilon}\,\frac{\mu_1 + \dots + \mu_n}{n}.\end{equation}

\noindent In particular, the normalisation $\int_X\beta_{\varepsilon}\wedge\gamma_{m,\,\varepsilon}^{n-1} = 1$ reads

\begin{equation}\label{eqn:normalisation-translation}\frac{1}{c_{\varepsilon}}\,\int\limits_X\lambda_1\dots\lambda_n\, \gamma_{m,\,\varepsilon}^n = \int\limits_X\frac{\mu_1 + \dots + \mu_n}{n}\, \gamma_{m,\,\varepsilon}^n =1.\end{equation}

\noindent Note that we also have

\begin{equation}\label{eqn:alpha_n-1_beta}\widetilde{\alpha}_{m,\,\varepsilon}^{n-1}\wedge\beta_{\varepsilon} = \frac{1}{n}\,(\Lambda_{\widetilde{\alpha}_{m,\,\varepsilon}}\beta_{\varepsilon})\,\widetilde{\alpha}_{m,\,\varepsilon}^n = \frac{1}{n}\,(\Lambda_{\widetilde{\alpha}_{m,\,\varepsilon}}\beta_{\varepsilon})\,\lambda_1\dots\lambda_n\,\gamma_{m,\,\varepsilon}^n.\end{equation}

 Putting all of the above together, we get\!\!:

\begin{eqnarray}\nonumber & & \bigg(\int\limits_X\widetilde{\alpha}_{m,\,\varepsilon}\wedge\gamma_{m,\,\varepsilon}^{n-1}\bigg)\cdot\bigg(\int\limits_X\widetilde{\alpha}_{m,\,\varepsilon}^{n-1}\wedge\beta_{\varepsilon}\bigg)\\
\nonumber & = & \bigg(\int\limits_X\frac{1}{n}\,(\Lambda_{\gamma_{m,\,\varepsilon}}\widetilde{\alpha}_{m,\,\varepsilon})\,\gamma_{m,\,\varepsilon}^n\bigg)\cdot\bigg(\int\limits_X \frac{1}{n}\,(\Lambda_{\widetilde{\alpha}_{m,\,\varepsilon}}\beta_{\varepsilon})\,\lambda_1\dots\lambda_n\,\gamma_{m,\,\varepsilon}^n\bigg)\\
\nonumber & \stackrel{(a)}{\geq} & \frac{1}{n^2}\,\bigg(\int\limits_X\bigg[(\Lambda_{\gamma_{m,\,\varepsilon}}\widetilde{\alpha}_{m,\,\varepsilon})\,(\Lambda_{\widetilde{\alpha}_{m,\,\varepsilon}}\beta_{\varepsilon})\bigg]^{\frac{1}{2}}\,(\lambda_1\dots\lambda_n)^{\frac{1}{2}}\,\gamma_{m,\,\varepsilon}^n\bigg)^2\\
\nonumber & \stackrel{(b)}{\geq} & \frac{1}{n^2}\,\bigg(\int\limits_X(\Lambda_{\gamma_{m,\,\varepsilon}}\beta_{\varepsilon})^{\frac{1}{2}}\,(\lambda_1\dots\lambda_n)^{\frac{1}{2}}\,\gamma_{m,\,\varepsilon}^n\bigg)^2  \stackrel{(c)}{=} \frac{1}{n^2}\,\bigg(\int\limits_X\frac{\sqrt{n}}{\sqrt{c_{\varepsilon}}}\,\lambda_1\dots\lambda_n\,\gamma_{m,\,\varepsilon}^n\bigg)^2\\
\nonumber  & = & \frac{1}{n\,c_{\varepsilon}}\,\bigg(\int\limits_X\widetilde{\alpha}_{m,\,\varepsilon}^n\bigg)^2 \stackrel{(d)}{=} \frac{1}{n\,c_{\varepsilon}}\,\bigg(\int\limits_Xc_{\varepsilon}\,\beta_{\varepsilon}\wedge\gamma_{m,\,\varepsilon}^{n-1}\bigg)^2 \stackrel{(e)}{=} \frac{c_{\varepsilon}}{n}.\end{eqnarray}

\noindent This proves (\ref{eqn:crucial-estimate}). Inequality $(a)$ is an application of the Cauchy-Schwarz inequality, inequality $(b)$ has followed from (\ref{eqn:trace-prod-ineq}), identity $(c)$ has followed from (\ref{eqn:eqn-translation}), identity $(d)$ has followed from $\widetilde{\alpha}_{m,\,\varepsilon}^n = c_{\varepsilon}\,\beta_{\varepsilon}\wedge\gamma_{m,\,\varepsilon}^{n-1}$ (which is nothing but the Monge-Amp\`ere equation $(\star)_{m,\,\varepsilon}$), while identity $(e)$ has followed from the normalisation $\int_X\beta_{\varepsilon}\wedge\gamma_{m,\,\varepsilon}^{n-1} = 1$ (cf. (\ref {eqn:contrad-ineq})). The proof of Lemma \ref{Lem:crucial-estimate} is complete.  \hfill $\Box$

\vspace{3ex}

\noindent {\it End of proof of Theorem \ref{The:main}.} Now, $\widetilde{\alpha}_{m,\,\varepsilon} = \alpha_{\varepsilon} + i\partial\bar\partial u_{m,\,\varepsilon}$ and $\partial\bar\partial\gamma_{m,\,\varepsilon}^{n-1}=0$, so

\begin{equation}\label{eqn:1st-term-ubound}\int\limits_X\widetilde{\alpha}_{m,\,\varepsilon}\wedge\gamma_{m,\,\varepsilon}^{n-1} =  \int\limits_X\alpha_{\varepsilon}\wedge\gamma_{m,\,\varepsilon}^{n-1}\longrightarrow \int\limits_X\alpha_{\varepsilon}\wedge\Gamma_{\infty,\,\varepsilon}\leq 1  \hspace{2ex} \mbox{for all}\hspace{1ex} \varepsilon>0,\end{equation}

\noindent where the above arrow stands for convergence as $m\rightarrow +\infty$ and the last inequality is nothing but (\ref{eqn:alpha-epsilon-int-bound}) (which, recall, is a consequence of the assumption that no K\"ahler current exists in $\{\alpha - \beta\}$ --- an assumption that we are going to contradict). On the other hand, the second factor on the l.h.s. of (\ref{eqn:crucial-estimate}) is precisely $M_{\varepsilon}$ defined in (\ref{eqn:M_0-def}), so in particular it is independent of $m$. Fixing any $\varepsilon>0$, taking limits as $m\rightarrow +\infty$ in (\ref{eqn:crucial-estimate}) and using (\ref{eqn:1st-term-ubound}), we get

\begin{equation}\label{eqn:m-limit-result}M_{\varepsilon}\geq\frac{c_{\varepsilon}}{n}   \hspace{3ex} \mbox{for every}\hspace{1ex}\varepsilon>0.\end{equation}

\noindent Taking now limits as $\varepsilon\downarrow 0$ and using (\ref{eqn:M_0-def}) and (\ref{eqn:c_0-def}), we get

$$M_0\geq\frac{c_0}{n}, \hspace{3ex} \mbox{i.e.} \hspace{2ex} \{\alpha\}^{n-1}.\,\{\beta\}\geq\frac{\{\alpha\}^n}{n}.$$

\noindent The last identity means that $\{\alpha\}^n - n\,\{\alpha\}^{n-1}.\,\{\beta\}\leq 0$ which is impossible if we suppose that $\{\alpha\}^n - n\,\{\alpha\}^{n-1}.\,\{\beta\}> 0$. This is the desired contradiction proving the existence of a K\"ahler current in the class $\{\alpha - \beta\}$ under the assumption $\{\alpha\}^n - n\,\{\alpha\}^{n-1}.\,\{\beta\}> 0$.  \hfill $\Box$

\vspace{3ex}

\noindent {\bf References} \\

\noindent [BDPP13]\, S. Boucksom, J.-P. Demailly, M. Paun, T. Peternell --- {\it The Pseudo-effective Cone of a Compact K\"ahler Manifold and Varieties of Negative Kodaira Dimension} --- J. Alg. Geom. {\bf 22} (2013) 201-248.

\vspace{1ex}

\noindent [Chi13]\, I. Chiose --- {The K\"ahler Rank of Compact Complex Manifolds} --- arXiv e-print CV 1308.2043v1

\vspace{1ex}

\noindent [Lam99]\, A. Lamari --- {Courants k\"ahl\'eriens et surfaces compactes} --- Ann. Inst. Fourier, Grenoble, {\bf 49}, 1 (1999), 263-285.

\vspace{1ex}

\noindent [Pop13]\, D. Popovici --- {\it Aeppli Cohomology Classes Associated with Gauduchon Metrics on Compact Complex Manifolds} --- arXiv e-print DG 1310.3685v1, to appear in Bull. Soc. Math. France.

\vspace{1ex}

\noindent [TW10]\, V. Tosatti, B. Weinkove --- {\it The Complex Monge-Amp\`ere Equation on Compact Hermitian Manifolds} --- J. Amer. Math. Soc. {\bf 23} (2010), no. 4, 1187-1195.

\vspace{1ex}

\noindent [Xia13]\, J. Xiao --- {\it Weak Transcendental Holomorphic Morse Inequalities on Compact K\"ahler Manifolds} --- arXiv e-print CV 1308.2878v1, to appear in Ann. Inst. Fourier.

\vspace{3ex}

\noindent Institut de Math\'ematiques de Toulouse, Universit\'e Paul Sabatier,

\noindent 118 route de Narbonne, 31062 Toulouse, France

\noindent Email\!: popovici@math.univ-toulouse.fr

\end{document}